\documentclass[12pt]{article}

\usepackage{amsmath}
\usepackage{amssymb}

\newcommand{\ci}{\mathrm{i}}
\newcommand{\dd}{\mathrm{d}}
\newcommand{\ee}{\mathrm{e}}

\usepackage{color}

\usepackage[hidelinks]{hyperref}

 \usepackage[margin=1in]{geometry} 

\usepackage{graphicx}
\usepackage{epstopdf}

\usepackage{caption}
\usepackage{subcaption}
\usepackage{float}

\def\XXint#1#2#3{{\setbox0=\hbox{$#1{#2#3}{\int}$}
     \vcenter{\hbox{$#2#3$}}\kern-.5\wd0}}

\title{Thermal contact resistance in two-phase media \\with arbitrarily shaped adiabatic fractures}
\author{E. Luca$^{1}$, J. S. Marshall$^{2}$ \& A. Chatziafratis$^{3,4}$}
\date{}

\begin{document}

\maketitle

\renewcommand{\thefootnote}{}
\footnotetext{{\it Date this work was submitted for peer-review: 12 June 2026}}

\begin{center}
$^1$Department of Mathematics and Statistics \\
University of Cyprus \\
Nicosia, Cyprus
\end{center}
\begin{center}
$^2$Department of Mathematics \\
Imperial College London \\
180 Queen's Gate \\
London SW7 2AZ, UK
\end{center}
\begin{center}
$^3$Waseda Research Institute for Science and Engineering \\
Waseda University \\
Tokyo, Japan
\end{center}
\begin{center}
$^4$Division of Applied Analysis \& PDE \\
Department of Mathematics \\
National \& Kapodistrian University of Athens \\
Athens, Greece
\end{center}
\begin{center}
$^{*}$Corresponding author: \tt{louca.elena@ucy.ac.cy}
\end{center}
\vspace{0.5cm}

\begin{abstract}
\noindent
This paper investigates thermal contact resistance in two-dimensional, periodic, two-phase media containing adiabatic fractures of arbitrary geometry. While traditional analytical models are often restricted to simplified contact shapes, such as flat or circular arcs, the current work employs the Unified Transform Method (Fokas, {\em PRSA}, {\bf 453}, 1997) to develop a generalized framework for calculating temperature fields and effective resistance across a broad class of geometries. The proposed method provides quasi-analytical solutions for the corresponding temperature fields and is validated against existing exact and approximate solutions for specific geometries. Its utility is further illustrated through several configurations featuring general fracture shapes. Our results provide a robust tool for predicting thermal performance in composite materials where interfacial roughness deviates from idealized geometries.
\end{abstract}


\newpage

\section{Introduction}

Heat conduction across interfaces between solid bodies is a classical problem in heat transfer, with a wide range of applications in engineering systems such as electronic packaging, energy devices, and composite materials \cite{DundursPanek:1976,MuzychkaYovanovich:2023,YovanovichMarotta:2003,Yovanovich:2005}. When two solids are brought into contact, the interface typically exhibits a temperature discontinuity due to imperfect thermal coupling. This phenomenon, commonly referred to as {\it thermal contact resistance} (or its inverse, {\it thermal contact conductance}) \cite{Razavietal:2016,YovanovichMarotta:2003}, is primarily caused by the presence of surface irregularities that limit the actual area of contact. This resistance depends not only on the intrinsic material properties but also on surface roughness, contact pressure, and the distribution of discrete contact spots, which typically occupy only a small fraction of the nominal contact area. Understanding and accurately modelling interfacial thermal resistance is therefore essential for predicting the effective thermal behaviour of heterogeneous media.

In 1969, Cooper {\em et al.}~\cite{CooperMikicYovanovich:1969} published a foundational study on contact resistance, deriving an estimate by solving for the axisymmetric temperature distribution within a cylindrical domain. Their model utilized a base characterized by an ``almost isothermal" contact zone surrounded by an adiabatic annulus, while maintaining an adiabatic outer boundary and uniform heat flux at the far field. A common limitation of the analysis in \cite{CooperMikicYovanovich:1969} and subsequent research is the assumption of a planar contact interface with periodic boundary conditions, which fails to account for the true topography of physical surface roughness. Addressing this, Das \& Sadhal \cite{DasSadhal:1992}  investigated a Cartesian-geometry problem involving flat contacts between semi-infinite media with differing thermal conductivities. Their model incorporated sparsely distributed, circular-arc gaps of varying protrusion angles, filled with a third material of finite conductivity. In a follow-up study \cite{DasSadhal:1998}, they extended this single-gap solution to a symmetric multi-gap configuration to account for interfacial interactions. However, to maintain analytical tractability, they assumed the gaps were sufficiently thin to simplify both the internal temperature profiles and the coupling between adjacent features.

Recent literature has highlighted a profound mathematical analogy between thermal contact resistance and the fluid dynamics of superhydrophobic surfaces \cite{HodesKirkCrowdy:2018,Mayeretal:2019,MiyoshiKirkHodes:2025}. Specifically, the problem of heat conduction through an isothermal contact flanked by adiabatic gaps is equivalent to the linear-shear flow of a liquid over a superhydrophobic surface in the Cassie (unwetted) state. In the hydrodynamic context, no-slip and shear-free boundary conditions are applied at the liquid-solid and liquid-gas interfaces (menisci), respectively. When appropriately non-dimensionalised, it has been shown that the dimensionless thermal contact resistance and the dimensionless spreading resistance are identical to the dimensionless apparent hydrodynamic slip length. This equivalence is particularly valuable because it allows for the adaptation of advanced results from the study of ``bubble mattresses" (surfaces where the liquid flows over curved menisci) \cite{Crowdy:2016} to thermal problems involving non-planar, curved adiabatic boundaries. For instance, the studies \cite{HodesKirkCrowdy:2018,Mayeretal:2019,MiyoshiKirkHodes:2025} demonstrated that recent developments accounting for meniscus curvature can be immediately applied to provide explicit formulae for contact resistance that capture both boundary curvature and contact distribution effects. 

\begin{figure}
\centering
\includegraphics[scale=1.2]{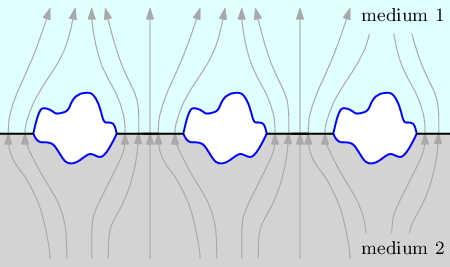}
\caption{Thermal contact resistance in two dimensions, in two-phase media with a periodic array of adiabatic fractures of arbitrary shape.}
\label{fig0}
\end{figure}

Recently, Jain \cite{Jain:2026} developed an exact theoretical framework for thermal contact resistance in the presence of multiple two-dimensional interfacial gaps with arbitrary, non-circular geometries. By modelling the gaps as regions of zero thermal conductivity within a semi-infinite flux-channel geometry, the method reformulates the problem as heat conduction in a heterogeneous medium with spatially varying conductivity, enabling an eigenfunction series representation for the temperature field. Unlike earlier asymptotic or perturbative approaches restricted to small protrusion angles or extreme contact fractions, the formulation in \cite{Jain:2026} remains applicable across a broad parameter range and can accommodate discrete profilometer data describing realistic surface roughness. The study also demonstrated the treatment of multiple irregular gaps, including triangular, sinusoidal, and experimentally measured roughness profiles, thereby significantly extending the class of geometries accessible to analytical treatment.

Despite their success, most existing models rely on the simplifying assumption that the contacting surfaces are nominally flat and that the interface geometry can be represented by regular or weakly perturbed configurations. Recent work by Jain \cite{Jain:2026} significantly broadened this framework by developing an eigenfunction-series formulation capable of treating multiple irregular interfacial gaps of arbitrary shape within a semi-infinite flux-channel geometry. In many practical situations, particularly in heterogeneous or composite media, interface geometries may be highly irregular. The presence of fractures, voids, or non-contact regions of arbitrary shape introduces significant geometric complexity, resulting in nontrivial temperature fields and altered heat flux pathways. These features can strongly influence the effective thermal resistance, yet they are not adequately captured by classical theories based on planar or axisymmetric geometries.

In this paper, we analyse thermal contact resistance in two dimensions, in two-phase media containing a periodic array of arbitrarily shaped adiabatic fractures (Figure \ref{fig0}). The fractures represent non-contact regions that interrupt heat flow across the interface, thereby modifying the effective resistance. The problem is analysed by implementing a suitable extension of the celebrated Unified Transform Method (UTM), also known as the Fokas method, which was pioneered by A. S. Fokas \cite{Fokas:1997} in the late 1990's as a general transform approach for initial-boundary value problems and subsequently developed in a variety of directions by Fokas and a large community of mathematicians \cite{fokasbook}. Although at first the method was chiefly introduced for integrable nonlinear dispersive PDE, relying on Riemann--Hilbert machinery and the Lax-pairs formalism, it soon also gave rise to a novel, effective approach for analysing and solving very general problems for linear evolution \cite{DeconinckTrogdonVasan:2014,Fokas:2002,Fokas:2009} as well as elliptic equations \cite{Ashton:2012,Fokas:2001}. The method has since been established as an essential generalisation of classical transforms (Fourier, Laplace) and solution methods (Green's functions, Wiener--Hopf technique), and may be viewed as a synthesis, rather than a separation, of variables \cite{FokasSpence:2012}: it formulates the problem in spectral space by combining the governing equation, boundary conditions and analyticity properties of suitable spectral functions, with the global relation providing the central equation coupling the transforms of known and unknown boundary values. Regarding elliptic problems, in particular, the UTM has been developed extensively, in the last three decades, for two-dimensional boundary value problems. Following Fokas's seminal work on linear PDE in convex polygons, substantial progress has been made pertaining to miscellaneous problems for Laplace's equation on polygonal and circular domains, curvilinear polygons, and more general planar geometries; see, for example, \cite{Colbrook:2020,ColbrookFlyerFornberg:2018,ColbrookFokasHash2019,Crowdy:2015,Crowdy:2015b,FokasKapaev:2003,HLLL:2024,HLLL:2025,LucaMarshallKaramanis:2018}. In light of the aforementioned advances, the UTM is deemed particularly well-suited to the present setting, ultimately providing a natural semi-analytical framework within which one may proceed efficiently. For the sake of completeness and the benefit of our readership, we conclude by mentioning indicatively also some recent works relating to dynamical configurations for variations of the traditional heat equation \cite{CACF:2024,CFA:2024,CLW:2026,CMKFA:2025} and to a time-dependent interface problem, posed on two adjacent half-spaces, for a non-local model describing water-waves propagation \cite{ChatziafratisOzawa:2026}.


By employing the UTM, we derive quasi-analytical expressions for the temperature fields in both media, accommodating complex fracture geometries. Having derived the solutions using the UTM, we report computed values for the increase in contact resistance for various fracture shapes, as well as temperature on the fracture boundaries. This approach allows a quantitative assessment of how interface geometry, fracture shape, and boundary conditions jointly influence thermal contact resistance. It is important to emphasise that using our approach, one can obtain an effective measure for the contact resistance without computing the full temperature field. This work provides a framework for bridging the gap between idealized contact models and realistic interface geometries, offering new insights into how fracture shape and distribution influence thermal transport in heterogeneous media. Our results are validated against analytical or approximate solutions in special cases.

\section{Problem formulation}

\begin{figure}
\centering
\includegraphics[scale=0.9]{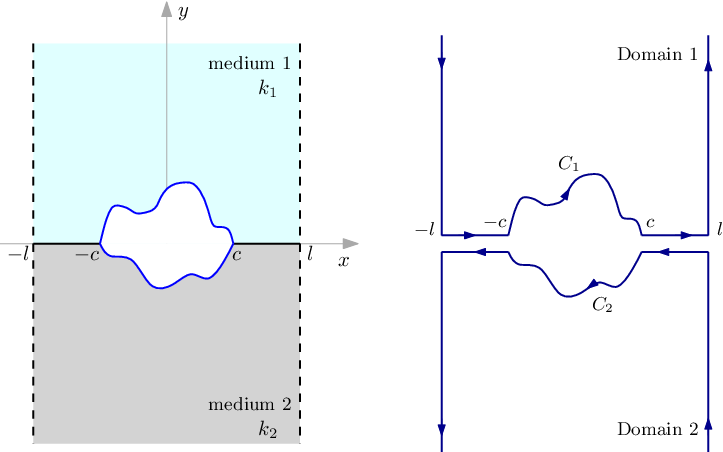}
\caption{Problem configuration for steady heat conduction through a two-phase medium with conductivities $k_1$ and $k_2$ in the upper and lower medium, respectively. A period window with an adiabatic fracture of arbitrary shape (left) and domain decomposition (right).}
\label{fig1}
\end{figure}

Under steady state conditions, the temperature field $T_{j}(x,y)$, $j=1,2$ in each medium (subscripts $1$ and $2$ refer to the upper and lower medium, respectively - see Fig. \ref{fig1}) satisfies Laplace's equation
\begin{equation}
\nabla^{2} T_{j}(x,y)=0.
\end{equation}
It is natural to introduce the usual complex variable $z=x+\ci y$ and complex potential function
\begin{equation}
w_{j}(z)=\chi_{j}(x,y)+\ci \, T_{j}(x,y), \quad j=1,2,
\end{equation}
in each medium, where $\chi_{j}(x,y)$ is the harmonic conjugate to $T_{j}(x,y)$. The far-field conditions are
\begin{equation}
\begin{split}\label{far-field conditions}
\frac{\partial T_{1}}{\partial y} \rightarrow -1, \quad \text{as } y \rightarrow +\infty, \\
\frac{\partial T_{2}}{\partial y} \rightarrow -r, \quad \text{as } y \rightarrow -\infty,
\end{split}
\end{equation}
where $r=k_{1}/k_{2}$ is the ratio of thermal conductivities. Note that, from \eqref{far-field conditions}, the heat flow rates at $y \rightarrow \pm \infty$ are consistent. Next, we write
\begin{equation}
\begin{split}\label{expansions}
w_{1}(z)&=h_1(z)+f_{1}(z), \qquad h_1(z)=-z+b_{1}, \\
w_{2}(z)&=h_2(z)+f_{2}(z), \qquad h_2(z)=-r z+b_{2},
\end{split}
\end{equation}
for some unknown constants $b_{1}, b_{2}$. The functions $f_{1}(z)$, $f_{2}(z)$ are analytic in their domains of definition and vanish as $y \rightarrow \infty$ and $y \rightarrow - \infty$ respectively. It can be shown that $h_j(z), j=1,2$, satisfy the far-field conditions \eqref{far-field conditions}.

{\bf Boundary and interface conditions:} The temperature field $T_{j}(x,y), j=1,2$ is $2l$-periodic in $x$ which implies that the complex potential satisfies the {\it quasi-periodicity} condition
\begin{equation}\label{quasi-period}
\begin{split}
w_{1}(z+2l)&=w_{1}(z)-2l, \\
w_{2}(z+2l)&=w_{2}(z)-2rl.
\end{split}
\end{equation}

The fracture is {\it adiabatic} and this implies that
\begin{equation}
\frac{\partial T_{1}}{\partial n}=\frac{\partial T_{2}}{\partial n}=0.
\end{equation}
On use of the Cauchy--Riemann equations, these conditions are equivalent to
\begin{equation}\label{chi1chi2}
\chi_{1}=\chi_{2}= C,
\end{equation}
where $C$ is a real constant. Without loss of generality we set $C$ to be equal to zero. Expression \eqref{chi1chi2} can then be written in terms of the complex potentials as
\begin{equation}\label{fracture BC}
\text{Re}[w_{1}(z)]=\text{Re}[w_{2}(z)]=0.
\end{equation}

On the contact interfaces $x \in [-l,-c] \cup [c,l]$, we have
\begin{equation}
T_{1}=T_{2}, \qquad \frac{\partial T_{1}}{\partial y}=p \, \frac{\partial T_{2}}{\partial y},
\end{equation}
where $p=1/r$. On use of the Cauchy--Riemann equations these are equivalent to
\begin{equation}
T_{1}=T_{2}, \qquad \chi_{1} = p \,\chi_{2},
\end{equation}
which can be written in terms of the complex potential functions as
\begin{equation}\label{interface conditions}
\text{Im}[w_{1}(z)]=\text{Im}[w_{2}(z)], \qquad  \text{Re}[w_{1}(z)]=p \,\text{Re}[w_{2}(z)].
\end{equation}

\section{Domain splitting}

The problem domain is split into two sub-domains (Fig.~\ref{fig1}) and in each sub-domain the complex potential is defined in \eqref{expansions}. The aim is to compute the unknown correction functions $f_{1}(z)$, $f_{2}(z)$ and parameters $b_{1}, b_{2}$.

\subsection{Domain 1}

We introduce the spectral functions $\rho_{j}(k)$, $j=1,2,3,4,5$ which are defined by
\begin{equation}
\begin{split}
\rho_{1}(k)&=\int_{l}^{l+\ci \infty} {f_{1}(z) \ee^{-\ci k z} \dd z}, \qquad \rho_{2}(k)=\int_{-l+\ci \infty}^{-l} {f_{1}(z) \ee^{-\ci k z} \dd z}, \\ 
\rho_{3}(k)&=\int_{-l}^{-c} {f_{1}(z) \ee^{-\ci k z} \dd z}, ~~~~~~~~~ ~~\rho_{4}(k)= \int_{C_{1}} {f_{1}(z) \ee^{-\ci k z} \dd z}, \\
\rho_{5}(k)&=\int_{c}^{l} {f_{1}(z) \ee^{-\ci k z} \dd z},
\end{split}\label{rho spectral functions}
\end{equation}
where contour $C_{1}$ is traversed as shown in Fig.~\ref{fig1}. The global relation is given by 
\begin{equation}\label{GR1}
\sum_{j=1}^{5} {\rho_{j}(k)}=0, \quad \text{Re} [k] \leq 0.
\end{equation}

\subsection{Domain 2}

We introduce the spectral functions $\sigma_{j}(k)$, $j=1,2,3,4,5$ which are defined by
\begin{equation}
\begin{split}
\sigma_{1}(k)&=\int_{l-\ci \infty}^{l} {f_{2}(z) \ee^{-\ci k z} \dd z}, \qquad \sigma_{2}(k)=\int_{-l}^{-l-\ci \infty} {f_{2}(z) \ee^{-\ci k z} \dd z}, \\
\sigma_{3}(k)&=\int_{-c}^{-l} {f_{2}(z) \ee^{-\ci k z} \dd z}, ~~~~~~~~~~ \sigma_{4}(k)=\int_{C_{2}} {f_{2}(z) \ee^{-\ci k z} \dd z}, \\
\sigma_{5}(k)&=\int_{l}^{c} {f_{2}(z) \ee^{-\ci k z} \dd z},
\end{split}
\end{equation}
where contour $C_{2}$ is traversed as shown in Fig. \ref{fig1}. The global relation is given by 
\begin{equation}\label{GR2}
\sum_{j=1}^{5} {\sigma_{j}(k)}=0, \quad \text{Re} [k] \geq 0.
\end{equation}

\subsection{Fracture shapes}

If the fracture is {\it flat}, then integral representations for $f_{1}(z), f_{2}(z)$ can be written following Fokas \& Kapaev \cite{FokasKapaev:2003}. In this case, we have
\begin{equation}
\rho_{4}(k)=\int_{-c}^{c} {f_{1}(z) \ee^{-\ci k z} \dd z},
\end{equation}
and the correction complex potential $f_{1}(z)$ can be represented by \cite{FokasKapaev:2003}:
\begin{equation}
f_{1}(z)=\frac{1}{2 \pi} \left[ \int_{0}^{-\ci \infty}{\rho_{1}(k) \ee^{\ci kz} \dd k} +\int_{0}^{\ci \infty}{\rho_{2}(k) \ee^{\ci kz} \dd k} + \sum_{j=3}^{5} \int_{0}^{\infty}{{\rho_{j}(k)} \ee^{\ci kz} \dd k} \right].
\end{equation}
A similar representation can be written for $f_2(z)$.

If the fracture in domain 1 has a {\it circular} arc shape, then an integral representations for $f_{1}(z)$ can be written following Crowdy \cite{Crowdy:2015}. We have
\begin{equation}
\rho_{4}(k)= \int_{C_{1}} {f_{1}(z) \ee^{-\ci k z} \dd z},
\end{equation}
where $C_{1}$ is a boundary part of a circle with centre $z_{0}$ and radius $q$, and which is traversed as shown in Fig.~\ref{fig1}. We also introduce
\begin{equation}
\rho_{0}(k)=\frac{1}{q} \int_{C_{1}} {f_{1}(z)  \left[\frac{z-z_{0}}{q}\right]^{k} \dd z}.
\end{equation}
The correction complex potential $f_{1}(z)$ can be represented by \cite{Crowdy:2015}:
\begin{equation}
\begin{split}
f_{1}(z)&=\frac{1}{2 \pi}\left[ \int_{0}^{-\ci \infty} { \rho_{1}(k) \ee^{\ci k z} \dd k} + \int_{0}^{\ci \infty} {\rho_{2}(k) \ee^{\ci k z} \dd k} + \int_{0}^{\infty} { [\rho_{3}(k)+\rho_{5}(k)] \ee^{\ci k z} \dd k}\right], \\
&-\frac{1}{2\pi\ci} \left[
      \int_{L_{1}} {\frac{\rho_{0}(k)}{1-\ee^{2\pi \ci k}} \left[ \frac{q}{z-z_{0}} \right]^{k+1} \dd k}
      + \int_{L_{2}} {\rho_{0}(k) \left[ \frac{q}{z-z_{0}} \right]^{k+1} \dd k}
\right. \\
&\qquad\qquad\left.
      + \int_{L_{3}} {\frac{\rho_{0}(k) \ee^{2\pi \ci k}}{1-\ee^{2 \pi \ci k}} \left[ \frac{q}{z-z_{0}} \right]^{k+1} \dd k}
\right],
\end{split}
\end{equation}
where the set $\{L_{j}| j=1,2,3\}$ constitutes the fundamental contour for circular edges \cite{Crowdy:2015}.

If the fracture has an {\it arbitrary shape}, then $f_1(z), f_2(z)$ can be computed directly from the Cauchy integral formula. We also mention the recent work by Hulse {\em et al.} \cite{HLLL:2024,HLLL:2025} who provide a new transform-based approach for Laplace's equation in non-convex and other planar domains, specifically the so-called Lipschitz domains. This work may be extended to unbounded domains and, therefore, be used in the present study, to provide integral representations for $f_1(z), f_2(z)$.


\section{Boundary/interface conditions and spectral analysis}

In this section, we analyse the boundary and interface conditions to obtain relations between the spectral functions.

\subsection{Domain 1}
 
On substitution of \eqref{expansions} into \eqref{quasi-period}, we find that, on the left vertical boundary,
\begin{equation}\label{1 period bc}
f_{1}(z)=f_{1}(z+2l).
\end{equation}
We multiply \eqref{1 period bc} by $\ee^{-\ci kz}$ and integrate along the left vertical boundary,
\begin{equation}
\int_{-l+\ci \infty}^{-l} {f_{1}(z) \ee^{-\ci kz} \dd z} =\int_{-l+\ci \infty}^{-l} {f_{1}(z+2l) \ee^{-\ci kz} \dd z}.
\end{equation}
This can be written in terms of the spectral functions as
\begin{equation}\label{spectral upper vertical}
\rho_{2}(k)=-\ee^{2\ci kl} \rho_{1}(k). 
\end{equation}

On the upper side of the fracture ($C_1$), we have
\begin{equation}
\text{Re}[w_{1}(z)]=0,
\end{equation}
which, on substitution of \eqref{expansions}, becomes
\begin{equation}\label{1 fracture bc}
f_{1}(z)+\overline{f_{1}(z)}= - \left[h_1(z)+\overline{h_1(z)}\right].
\end{equation}
We multiply \eqref{1 fracture bc} by $\ee^{-\ci kz}$ and integrate along the upper part of the fracture $C_{1}$:
\begin{equation}
\int_{C_{1}} {f_{1}(z) \ee^{-\ci kz} \dd z} +\int_{C_{1}} {\overline{f_{1}(z)}\ee^{-\ci kz} \dd z}=D_{1}(k),
\end{equation}
where
\begin{equation}
D_{1}(k)= -\int_{C_{1}} {\left[h_1(z)+\overline{h_1(z)}\right] \ee^{-\ci kz} \dd z}.
\end{equation}
This can equivalently be written as
\begin{equation}\label{spectral upper fracture}
\rho_{4}(k) =-\int_{C_{1}} {\overline{f_{1}(z)}\ee^{-\ci kz} \dd z}+D_{1}(k).
\end{equation}

\subsection{Domain 2}

Similarly to the analysis above, we find that, on the left vertical boundary,
\begin{equation}\label{2 period bc}
f_{2}(z)=f_{2}(z+2l).
\end{equation}
We multiply \eqref{2 period bc} by $\ee^{-\ci kz}$ and integrate along the left vertical boundary, to obtain
\begin{equation}\label{spectral lower vertical}
\sigma_{2}(k)=-\ee^{2\ci kl} \sigma_{1}(k). 
\end{equation}

On the lower side of the fracture ($C_2$), we have
\begin{equation}
\text{Re}[w_{2}(z)]=0,
\end{equation}
which, on substitution of \eqref{expansions}, becomes
\begin{equation}\label{2 fracture bc}
f_{2}(z)+\overline{f_{2}(z)} = - \left[h_2(z)+\overline{h_2(z)}\right].
\end{equation}
We multiply \eqref{2 fracture bc} by $\ee^{-\ci kz}$ and integrate along the lower part of the fracture $C_{2}$, to find
\begin{equation}\label{spectral lower fracture}
\sigma_{4}(k) =-\int_{C_{2}} {\overline{f_{2}(z)} \ee^{-\ci kz} \dd z}+D_{2}(k),
\end{equation}
where
\begin{equation}
D_{2}(k)= - \int_{C_{2}} {\left[h_2(z)+\overline{h_2(z)}\right] \ee^{-\ci kz} \dd z}.
\end{equation}

\subsection{Interface}
On substitution of \eqref{expansions} into the first interface condition in \eqref{interface conditions}, we find
\begin{equation}\label{interface bc 1}
f_{1}(z)-\overline{f_{1}(z)} - f_{2}(z)+\overline{f_{2}(z)}=-h_1(z)+\overline{h_1(z)}+h_2(z)-\overline{h_2(z)}.
\end{equation}
We multiply \eqref{interface bc 1} by $\ee^{-\ci kz}$ and integrate along $x \in [-l,-c]$, to obtain
\begin{equation}\label{spectral interface 1L}
\rho_{3}(k)-\overline{\rho_{3}}(-k)+\sigma_{3}(k)-\overline{\sigma_{3}}(-k)=P(k),
\end{equation}
where
\begin{equation}
P(k)=\int_{-l}^{-c} {\left[-h_1(z)+\overline{h_1(z)}+h_2(z)-\overline{h_2(z)}\right] \ee^{-\ci kz} \dd z}.
\end{equation} 
Similarly, along $x \in [c,l]$, we find that
\begin{equation}\label{spectral interface 1R}
\rho_{5}(k)-\overline{\rho_{5}}(-k)+\sigma_{5}(k)-\overline{\sigma_{5}}(-k)=R(k),
\end{equation}
where
\begin{equation}
R(k)=\int_{c}^{l} {\left[-h_1(z)+\overline{h_1(z)}+h_2(z)-\overline{h_2(z)}\right] \ee^{-\ci kz} \dd z}.
\end{equation}

On substitution of \eqref{expansions} into the second interface condition in \eqref{interface conditions}, we find
\begin{equation}\label{interface bc 2}
f_{1}(z)+\overline{f_{1}(z)} -p \,\left[f_{2}(z) + \overline{f_{2}(z)}\right]= - h_1(z) -\overline{h_1(z)} + p \, \left[h_2(z)+\overline{h_2(z)}\right].
\end{equation}
Multiplication of \eqref{interface bc 2} by $\ee^{-\ci kz}$ and integration along $x \in [-l,-c]$ gives
\begin{equation}\label{spectral interface 2L}
\rho_{3}(k)+\overline{\rho_{3}}(-k)+ p \,\sigma_{3}(k) +p \,\overline{\sigma_{3}}(-k)=S(k),
\end{equation}
where
\begin{equation}
S(k)=\int_{-l}^{-c} {\left[-h_1(z)-\overline{h_1(z)}+p \,h_2(z) + p \,\overline{h_2(z)}\right] \ee^{-\ci kz} \dd z}.
\end{equation}
Along $x \in [c,l]$, we find that
\begin{equation}\label{spectral interface 2R}
\rho_{5}(k)+\overline{\rho_{5}}(-k)+p\, \sigma_{5}(k) + p\, \overline{\sigma_{5}}(-k)=T(k),
\end{equation}
where
\begin{equation}
T(k)=\int_{c}^{l} {\left[-h_1(z)-\overline{h_1(z)}+p\, h_2(z) + p\, \overline{h_2(z)}\right] \ee^{-\ci kz} \dd z}.
\end{equation}

\section{Solution scheme} 

In this section, we show how to exploit the analyticity properties of the spectral functions $\rho_1(k)$ and $\sigma_1(k)$ to identify special points in the spectral $k$-plane whereby information on a reduced set of unknown spectral functions is available. The spectral analysis of the boundary and interface conditions produced expressions \eqref{spectral upper vertical}, \eqref{spectral upper fracture}, \eqref{spectral lower vertical}, \eqref{spectral lower fracture}, \eqref{spectral interface 1L}, \eqref{spectral interface 1R}, \eqref{spectral interface 2L}, \eqref{spectral interface 2R}, to be analysed together with the global relations \eqref{GR1} and \eqref{GR2}.

\subsection{Conditions in the spectral $k$-plane}
On use of \eqref{spectral upper vertical} and \eqref{spectral upper fracture}, the global relation \eqref{GR1} can be written as
\begin{equation}
(1-\ee^{2\ci kl}) \rho_{1}(k)=W_{1}(k), \quad \text{Re}[k] \leq 0,
\end{equation}
where
\begin{equation}\label{W1 initial}
W_{1}(k)=-\rho_{3}(k)+\int_{C_{1}} {\overline{f_{1}(z)}\ee^{-\ci kz} \dd z}-D_{1}(k)-\rho_{5}(k).
\end{equation}
But $\rho_{1}(k)$ is analytic in the left half-plane which implies that $W_{1}(k)$ must vanish at zeros of $1-\ee^{2\ci kl}=0$ in the left half-plane, that is, at points in the set:
\begin{equation}
\Sigma_{1}=\left\{k_{n}= n\pi/l \bigl | n \in \mathbb{Z}, n \leq 0 \right\}.
\end{equation}
Therefore, we have the following set of conditions
\begin{equation}\label{conditions 1}
W_{1}(k)=0, \quad \text{for } k \in \Sigma_{1}.
\end{equation}

Similarly, on use of \eqref{spectral lower vertical} and \eqref{spectral lower fracture}, the global relation \eqref{GR2} can be written as
\begin{equation}
(1-\ee^{2\ci kl}) \sigma_{1}(k)=W_{2}(k), \quad \text{Re} [k] \geq 0,
\end{equation}
where
\begin{equation}\label{W2 initial}
W_{2}(k)=-\sigma_{3}(k)+\int_{C_{2}} {\overline{f_{2}(z)} \ee^{-\ci kz} \dd z}-D_{2}(k)-\sigma_{5}(k).
\end{equation}
But $\sigma_{1}(k)$ is analytic in the right half-plane which implies that $W_{2}(k)$ must vanish at zeros of $1-\ee^{2\ci kl}=0$ in the right half-plane, that is, at points in the set:
\begin{equation}
\Sigma_{2}=\left\{k_{n}= n\pi/l \bigl | n \in \mathbb{Z}, n \geq 0 \right\}.
\end{equation}
Therefore, we have the following set of conditions
\begin{equation}\label{conditions 2}
W_{2}(k)=0, \quad \text{for } k \in \Sigma_{2}.
\end{equation}

Using \eqref{spectral interface 1L}, \eqref{spectral interface 1R}, \eqref{spectral interface 2L}, \eqref{spectral interface 2R}, the spectral functions $\sigma_{3}(k)$ and $\sigma_{5}(k)$ can be expressed in terms of $\rho_{3}(k), \overline{\rho_{3}}(-k), \rho_{5}(k), \overline{\rho_{5}}(-k)$ via
\begin{equation}\label{tau3tau5}
\begin{split}
\sigma_{3}(k)&= - \frac{r+1}{2} \,\rho_{3}(k) - \frac{r-1}{2} \, \overline{\rho_{3}}(-k) + \frac{1}{2} \, [P(k)+r\,S(k)], \\
\sigma_{5}(k)&= - \frac{r+1}{2} \,\rho_{5}(k) - \frac{r-1}{2} \, \overline{\rho_{5}}(-k) + \frac{1}{2} \, [R(k) + r\,T(k)].
\end{split}
\end{equation}

\subsection{Function representation}

The unknown boundary data along the contact interfaces, as well as upper and lower sides of the fracture will be represented using Chebyshev expansions.

{\it Interface $x \in [-l,-c]$:}
This can be parametrised by
\begin{equation}
z_3(s)=-\frac{l+c}{2}+\left(\frac{l-c}{2}\right)s, \quad \text{for } s \in [-1,1].
\end{equation}
We write
\begin{equation}\label{sum1}
f_1(z_3(s))=\sum_{m=0}^{\infty} {a_{m} T_{m}(s)},
\end{equation}
for some set of coefficients $\{ a_{m} \in \mathbb{C} | m=0,1,2,\dots \}$ to be found, and where $T_{m}(s)=\cos(m \cos^{-1}s)$, $s \in [-1,1]$ is the $m$th Chebyshev polynomial. It follows that
\begin{equation}\label{representationrho3}
\rho_{3}(k)=\sum_{m=0}^{\infty} {a_{m} [A(k,m)]}, \quad A(k,m)=\int_{-1}^{1} {T_{m}(s) \ee^{-\ci k z_3(s)} z_3'(s) \dd s}.
\end{equation}

{\it Interface $x \in [c,l]$:}
This can be parametrised by
\begin{equation}
z_5(s)=\frac{l+c}{2}+\left(\frac{l-c}{2}\right)s, \quad \text{for } s \in [-1,1].
\end{equation}
We write
\begin{equation}\label{sum2}
f_1(z_5(s))=\sum_{m=0}^{\infty} {b_{m} T_{m}(s)},
\end{equation}
for some set of coefficients $\{ b_{m} \in \mathbb{C} | m=0,1,2,\dots \}$. It follows that
\begin{equation}\label{representationrho5}
\rho_{5}(k)=\sum_{m=0}^{\infty} {b_{m} [B(k,m)]}, \quad B(k,m)=\int_{-1}^{1} {T_{m}(s) \ee^{-\ci k z_5(s)} z_5'(s) \dd s}.
\end{equation}

{\it Fracture part $C_{1}$ (upper side):}
Assuming that this can be parametrised by $z_u(s)$, $s \in [-1,1]$, such that $z_u(-1)=-c$ and $z_u(1)=c$, we write
\begin{equation}\label{sum3}
f_1(z_u(s))=\sum_{m=0}^{\infty} {c_{m} T_{m}(s)},
\end{equation}
for some set of coefficients $\{ c_{m} \in \mathbb{C} | m=0,1,2,\dots \}$. It follows that
\begin{equation}\label{representationrho4}
\rho_{4}(k)=\sum_{m=0}^{\infty} {c_{m} [F(k,m)]}, \quad F(k,m)= \int_{-1}^{1} {T_{m}(s) \ee^{-\ci k z_u(s)} z'_u(s) \dd s}.
\end{equation}

{\it Fracture part $C_{2}$ (lower side):}
Similarly, assuming that this boundary can be parametrised by $z_l(s)$, $s \in [-1,1]$, such that $z_l(-1)=c$ and $z_l(1)=-c$, we write
\begin{equation}\label{sum4}
f_2(z_l(s))=\sum_{m=0}^{\infty} {d_{m} T_{m}(s)},
\end{equation}
for some set of coefficients $\{ d_{m} \in \mathbb{C} | m=0,1,2,\dots \}$. It follows that
\begin{equation}\label{representationsigma4}
\sigma_{4}(k)=\sum_{m=0}^{\infty} {d_{m} [G(k,m)]}, \quad G(k,m)= \int_{-1}^{1} {T_{m}(s) \ee^{-\ci k z_l(s)} z'_l(s) \dd s}.
\end{equation}

\subsection{Formulation of a linear system}

Using the representations \eqref{representationrho3}, \eqref{representationrho5}, \eqref{representationrho4} and \eqref{representationsigma4}, conditions \eqref{conditions 1} can be expressed as
\begin{equation}\label{linear 1}
\sum_{m=0}^{M_1} {a_{m} [A(k,m)]}+\sum_{m=0}^{M_2} {b_{m} [B(k,m)]}+\sum_{m=0}^{M_3} {\overline{c_{m}} [-F(k,m)]}+D_{1}(k)=0, \quad \text{for } k \in \Sigma_{1}.
\end{equation}

Similarly, conditions \eqref{conditions 2} can be expressed as
\begin{equation}
\begin{split}\label{linear 2}
&\sum_{m=0}^{M_{1}} {a_{m} \left[(r+1) \,A(k,m)\right]}+ \sum_{m=0}^{M_{1}} {\overline{a_{m}} \left[(r-1) \,A(k,m) \right]} \\
+&\sum_{m=0}^{M_{2}} {b_{m} \left[(r+1) \, B(k,m)\right]}+ \sum_{m=0}^{M_{2}} {\overline{b_{m}} \left[(r-1) \, B(k,m) \right]} \\
+&\sum_{m=0}^{M_{4}} {\overline{d_{m}} [2\,G(k,m)]} - 2\,D_{2}(k) + U(k)=0, \quad \text{for } k \in \Sigma_{2},
\end{split}
\end{equation}
where
\begin{equation}
\begin{split}
U(k)&= - P(k) - R(k) - r [S(k)+T(k)] \\
&= 2 \left[(r+1)\, b_{1} + (r-1) \,\overline{b_{1}} - 2 b_{2}\right]\left\{ \begin{array}{lll}
\dfrac{\sin(kl)-\sin(kc)}{k}, \quad &k \neq 0, \\ \\
l-c, \quad &k=0.
\end{array}\right.
\end{split}
\end{equation}
Note that the series representations \eqref{sum1}, \eqref{sum2}, \eqref{sum3} and \eqref{sum4} were truncated to include only terms up to $m=0,...,M_{i}$, $i=1,2,3,4$ respectively for suitable $M_{i}$. We have also used that $\overline{A}(-k,m)=A(k,m)$ and $\overline{B}(-k,m)=B(k,m)$ in \eqref{linear 2}.

There is an additional degree of freedom in the definition of temperature, which means we can specify it at a certain point in the problem domain. Without loss of generality, we assume that
\begin{equation}
T_1=\text{Im}[w_{1}(z)]=0, \quad \text{at } z=-c.
\end{equation}
Using \eqref{expansions} and \eqref{sum1}, this can be expressed as
\begin{equation}\label{linear 3}
\sum_{m=0}^{M_{1}} {a_{m} T_{m}(1)} -\sum_{m=0}^{M_{1}} {\overline{a_{m}} T_{m}(1)}+b_{1}-\overline{b_{1}}=0.
\end{equation}

We formulate a linear system for the unknown coefficients $\{ a_m, b_m, c_m, d_m \}$ and parameters $b_{1}, b_{2}$. The linear system comprises conditions \eqref{linear 1} and \eqref{linear 2} evaluated at points in the sets $\Sigma_1, \Sigma_2$ respectively, together with the condition \eqref{linear 3}. We found that the coefficients $\{ a_m, b_m, c_m, d_m \}$ decay quickly and, therefore, we choose the truncation parameter to be $M_i=M=16$, for $i=1,2,3,4$. Once the coefficients and parameters $b_1, b_2$ are found, $f_1(z)$ and $f_2(z)$ can be computed in the domain interior.

\section{Results}

In this section, we present solutions for different fracture shapes and fractions $c/l$. We define the normalised parameter $\lambda_{s}$ given by
\begin{equation}\label{normalised lambdas}
\lambda_{s}=\frac{\text{Im}[b_{2}-b_{1}]}{2c^{2}/l}.
\end{equation}
The parameters $b_1, b_2$ are found from the solution of the linear system \eqref{linear 1}, \eqref{linear 2} and \eqref{linear 3}. The presence of the adiabatic fracture induces a macroscopic temperature discontinuity across the interface. This jump is quantified by $\text{Im}[b_{2}-b_{1}]$, which defines an effective thermal resistance of the fractured interface. The quantity $\lambda_s$ may therefore be interpreted as an {\it effective thermal resistance measure} associated with the fracture, where geometric and interfacial details are condensed into a single parameter governing the macroscopic flow. In analogy with superhydrophobic surface theory \cite{LucaMarshallKaramanis:2018}, $\lambda_s$ plays a role analogous to an effective slip length.

\begin{figure}
\centering
\includegraphics[scale=0.38]{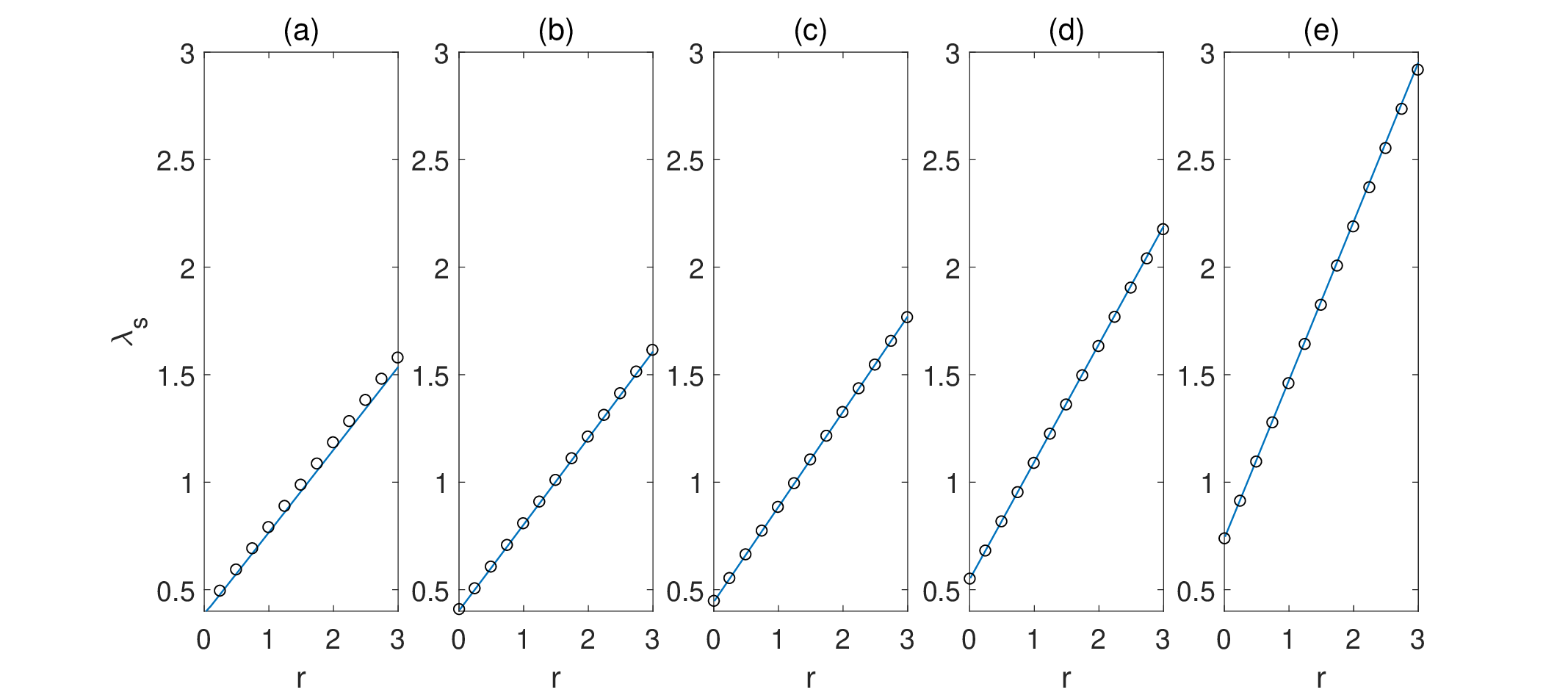}
\caption{Normalised parameter $\lambda_{s}$ given by \eqref{normalised lambdas} for a flat fracture $z_u(s)=cs=-z_l(s)$, $s\in[-1,1]$, as a function of the ratio of thermal conductivities $r=k_{1}/k_{2}$. Cases: (a) $c/l=0.1$, (b) $c/l=0.25$, (c) $c/l=0.5$, (d) $c/l=0.75$ and (e) $c/l=0.9$. Solid lines show our solutions and circles show the exact values found using conformal mapping.}
\label{fig:flat case}
\end{figure}

For a flat fracture, an exact solution can be obtained using conformal mapping techniques. As observed by \cite{HodesKirkCrowdy:2018}, \cite{Mayeretal:2019}, and \cite{MiyoshiKirkHodes:2025}, the present mathematical problem is equivalent to the problem of linear shear flow over a flat surface consisting of a periodic array of no-slip and shear-free strips. Consequently, the solution for a flat adiabatic fracture may be identified with the classical solution of Philip \cite{Philip:1972}, section 5, later rederived by Crowdy \cite{Crowdy:2011}, section VII. It is the latter form that we make use of here. This may be stated as follows. We introduce the conformal mapping from the unit disc
\begin{equation}\label{conformal mapping}
z(\zeta)=-\frac{\ci l}{\pi} \log \left[ \frac{(\zeta-\alpha)(\zeta-1/\alpha)}{(\zeta-\overline{\alpha}) (\zeta-1/\overline{\alpha})} \right],
\end{equation}
where $\alpha=\ci s$, $0<s<1$, with $s$ governing the fracture size. The upper-half unit $\zeta$-disc is mapped to domain 1, and the lower-half to domain 2. As $\zeta \rightarrow \alpha$ (which corresponds to $y \rightarrow \infty$), we have
\begin{equation}
z \sim - \frac{\ci l}{\pi} \log(\zeta - \alpha),
\end{equation}
while as $\zeta \rightarrow \overline{\alpha}$ (which corresponds to $y \rightarrow - \infty$), we have
\begin{equation}
z \sim \frac{\ci l}{\pi} \log(\zeta - \overline{\alpha}).
\end{equation}

It can be shown that the complex potentials $w_1(z)$ and $w_2(z)$ satisfying the far-field \eqref{far-field conditions}, boundary and interface conditions \eqref{quasi-period}, \eqref{fracture BC}, \eqref{interface conditions}, are given by
\begin{equation}\label{complex potentials exact}
\begin{split}
w_1(\zeta) &= \frac{\ci l}{\pi} \log \left[ \frac{(\zeta-\alpha)(\zeta-1/\overline{\alpha})}{(\zeta-\overline{\alpha}) (\zeta-1/\alpha)} \right], \\
w_2(\zeta) &= - r w_1(\zeta).
\end{split}
\end{equation}
Using the limiting behaviour of the conformal mapping \eqref{conformal mapping} as $\zeta \rightarrow \alpha$ and $\zeta \rightarrow \overline{\alpha}$, as well as the expressions for the complex potentials \eqref{complex potentials exact}, we find
\begin{equation}\label{b1 b2}
b_{1}=\frac{2 \ci l}{\pi} \log\left[\frac{\alpha^2 + 1}{\alpha^2 - 1}\right], \qquad b_{2}=\frac{2 \ci l r}{\pi} \log\left[\frac{\alpha^2 - 1}{\alpha^2 + 1}\right].
\end{equation}

Figure \ref{fig:flat case} shows the normalised parameter $\lambda_s$ for a flat fracture as a function of the ratio of thermal conductivities $r=k_1/k_2$, for different fractions $c/l$. Our solutions are verified against the exact solution for $\lambda_{s}$ given by \eqref{normalised lambdas} and \eqref{b1 b2}. The results demonstrate excellent agreement between the present UTM-based formulation and the exact conformal-mapping solution across the full range of conductivity ratios and fracture fractions considered. As expected, the effective resistance measure $\lambda_s$ increases as the fracture fraction $c/l$ increases, since larger adiabatic regions obstruct heat transport more strongly.

\begin{figure}
\centering
\begin{subfigure}{.35\textwidth}
  \centering
  \includegraphics[scale=0.8]{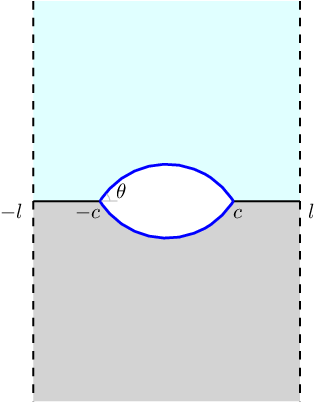}
\caption*{Fracture shape}
  \label{fig:thetap}
\end{subfigure}%
\begin{subfigure}{.7\textwidth}
  \centering
  \includegraphics[scale=0.6]{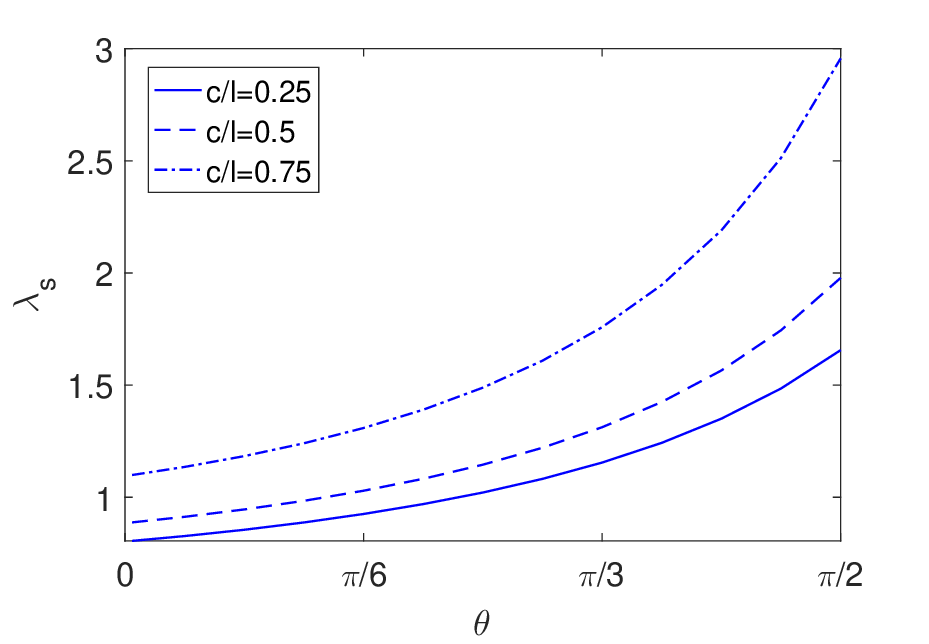}
  \label{fig:thetan}
\end{subfigure}
\caption{Normalised parameter $\lambda_{s}$ given by \eqref{normalised lambdas} for an up-down symmetric circular fracture as a function of angle $\theta$, for $r=k_{1}/k_{2}=1$ and $c/l=0.25, 0.5, 0.75$. The curve in domain 1 is described by $C_{1}: z_u(s)=z_{1}(s)=z_{0} + \ci r_0 \ee^{{\rm i}\theta s}$, $s \in [-1,1]$ ($(z_{0},r_0)$ is the (centre,radius) of a circular disc), where $z_0 \in \ci \mathbb{R}$. The curve in domain 2 is described by $C_{2}: z_l(s)=z_{2}(s)= - z_1(s)$, $s \in [-1,1]$.}
\label{fig:test2}
\end{figure}

Figure \ref{fig:test2} shows the normalised parameter $\lambda_{s}$ as a function of $\theta$ for an up-down symmetric circular fracture. Figure \ref{fig:test2} shows that the effective resistance measure depends strongly on the fracture curvature through the angle $\theta$. Increasing $\theta$ corresponds to increasingly protruding fracture geometries, which alter the local heat-flux pathways and enhance thermal obstruction. The effect becomes more pronounced as the fracture fraction $c/l$ increases, indicating stronger geometric interaction between neighbouring fractures.

Additional validation of the present formulation may also be inferred from the close mathematical analogy between thermal contact resistance and shear flow over superhydrophobic surfaces discussed by Hodes {\em et al.} \cite{HodesKirkCrowdy:2018}. In our earlier work on superhydrophobic surfaces \cite{LucaMarshallKaramanis:2018}, the same UTM-based methodology was applied to periodic slip/no-slip configurations and shown to be in excellent agreement with existing asymptotic and approximate theories. Owing to the equivalence between the dimensionless thermal resistance problem and the corresponding hydrodynamic slip problem, the present results for up-down symmetric fractures may therefore also be viewed as indirect validation against the approximate formulae developed for superhydrophobic surfaces by Crowdy \cite{Crowdy:2016}, Teo \& Khoo \cite{TeoKhoo:2010}, and Schnitzer \cite{Schnitzer:2017}. This correspondence further supports the accuracy and robustness of the present UTM framework across a broad range of interface geometries.

\begin{figure}
\centering
\begin{subfigure}{.35\textwidth}
  \centering
  \includegraphics[scale=0.8]{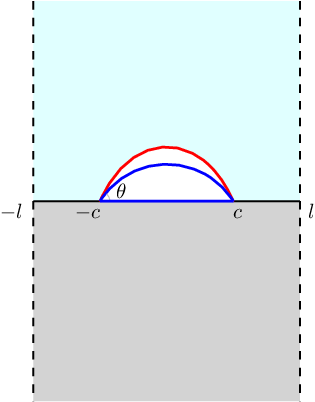}
\caption*{Fracture shape}
  \label{fig:thetap}
\end{subfigure}%
\begin{subfigure}{.7\textwidth}
  \centering
  \includegraphics[scale=0.52]{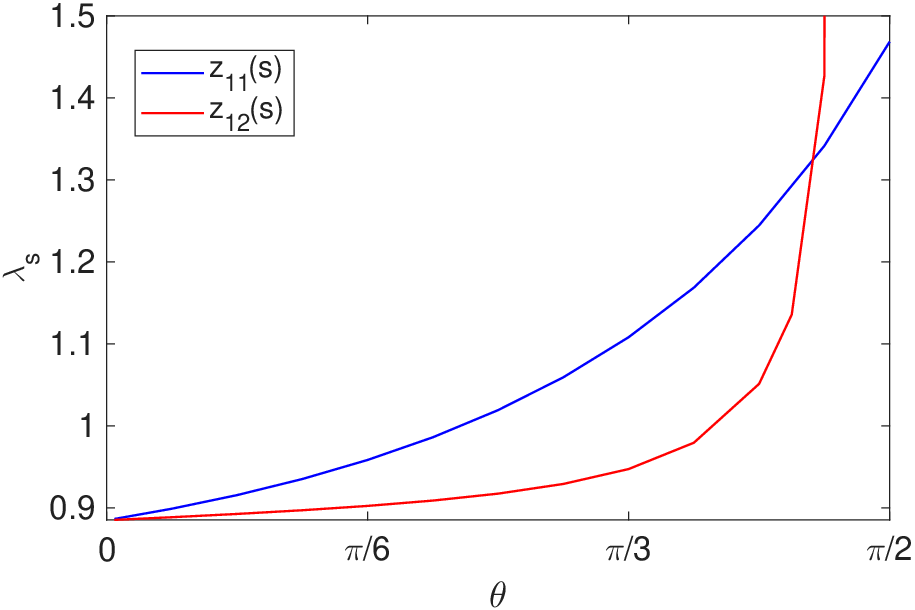}
  \label{fig:thetan}
\end{subfigure}
\caption{Normalised parameter $\lambda_{s}$ given by \eqref{normalised lambdas} for an asymmetric fracture as a function of angle $\theta$, for $r=k_{1}/k_{2}=1$ and $c/l=0.5$. The curve in blue (domain 1) is described by $C_{1}: z_u(s)=z_{11}(s)=z_{0} + \ci r_0 \ee^{{\rm i}\theta s}$, $s \in [-1,1]$ ($(z_{0},r_0)$ is the (centre,radius) of a circular disc). The curve in red is described by $C_{1}: z_u(s)=z_{12}(s)=- cs + \ci \beta \cos(\pi s/2)$, $s \in [-1,1]$, where $\beta=(2/\pi) \tan\theta$ is used to control angle $\theta$. The flat boundary in domain 2 is given by $z_l(s) = -cs$, $s\in[-1,1]$.}
\label{fig:test1}
\end{figure}

Figures \ref{fig:test1} shows the normalised parameter $\lambda_{s}$ as a function of $\theta$ for an asymmetric fracture. The results in Figure \ref{fig:test1} illustrate the sensitivity of thermal resistance to geometric asymmetry. Although the fractures possess comparable characteristic widths, the distinct upper-boundary geometries lead to noticeably different values of $\lambda_s$, demonstrating that the resistance depends not only on fracture size but also on detailed interfacial shape


%
%

\begin{figure}
\centering
\begin{subfigure}{.35\textwidth}
  \centering
  \includegraphics[scale=0.8]{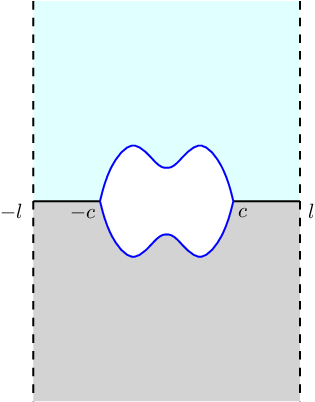}
\caption*{Fracture shape}
  \label{fig:thetap}
\end{subfigure}%
\begin{subfigure}{.7\textwidth}
  \centering
  \includegraphics[scale=0.6]{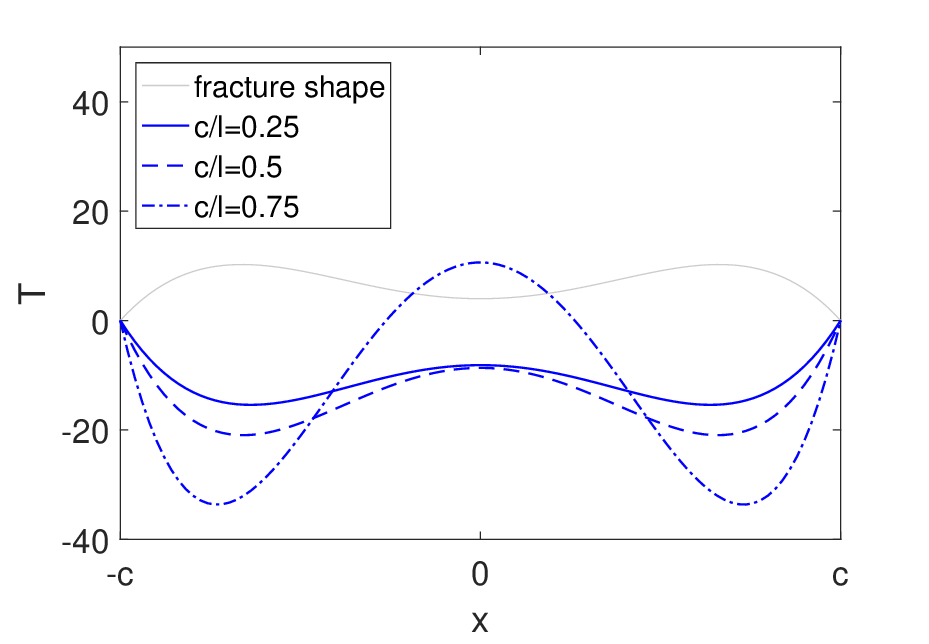}
  \label{fig:thetan}
\end{subfigure}
\caption{Temperature along the fracture boundary in domain 1 as a function of $x$, for $r=k_{1}/k_{2}=1$ and $c/l=0.25, 0.5, 0.75$. The curve in domain 1 is described by $C_{1}: z_u(s)=z_{1}(s)=x_{1}(s)+{\rm i} y_{1}(s)$, where $x_{1}(s)=-cs$, $y_{1}(s)=-(cs+2)(cs+1)(cs-1)(cs-2)+8$, $s \in [-1,1]$. The curve in domain 2 is described by $C_{2}: z_{l}(s)=-z_{u}(s)$, $s \in [-1,1]$. Constant $c$ is found using Newton's method.}
\label{fig:test3}
\end{figure}

Figure \ref{fig:test3} shows the temperature distribution along the fracture boundary for different fracture fractions $c/l$. We observe that the extrema occur near regions of highest boundary curvature, where the local redirection of heat flux is strongest. As the fracture fraction increases, the temperature variation along the boundary becomes increasingly non-uniform, indicating a stronger perturbation of the surrounding heat flux field by the enlarged adiabatic region. In particular, the larger values of $c/l$ produce sharper temperature gradients near the narrow neck regions of the fracture, reflecting enhanced flux crowding around these geometrically constrained portions of the interface. The symmetry of the temperature profiles about $x=0$ is consistent with the up-down and left-right symmetry of the fracture geometry.

Overall, the results demonstrate that both fracture geometry and fracture fraction $c/l$ have a significant influence on the effective thermal resistance. The UTM framework accurately reproduces known exact solutions while remaining flexible to treat complex symmetric and asymmetric fracture geometries beyond the reach of existing methods.

\section{Discussion}

The results presented in this work demonstrate that the UTM provides a robust framework for evaluating thermal contact resistance in media with complex interfacial geometries. A significant advantage of the current approach is its inherent flexibility; while the present analysis focuses on a single periodic fracture, the method can be naturally generalized to account for multiple fractures within a single periodic window. This extension is particularly relevant for modelling surfaces with hierarchical roughness or multi-scale contact distributions where a single-fracture approximation may be insufficient.

Furthermore, the UTM offers a distinct computational and analytical advantage over classical conformal mapping techniques. While conformal maps are powerful for specific two-dimensional geometries, they become increasingly difficult to implement, or altogether unavailable, when the interface morphology involves highly irregular shapes or multiple connected components. The current approach bypasses these geometric constraints by operating in the transform space, providing a unified path to the temperature field. 

\begin{figure}
\centering
\includegraphics[scale=0.8]{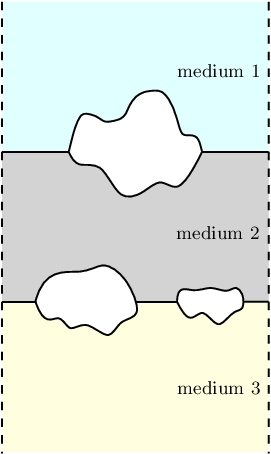}
\caption{Periodic array of adiabatic fractures of arbitrary shape in a three-phase media.}
\label{fig_3phase}
\end{figure}

Beyond the thermal domain, this methodology can be readily applied to other two-phase harmonic problems governed by Laplace's equation. Potential applications include predicting the effective electrical conductivity of composites with insulating inclusions or determining the longitudinal slip length in microfluidic channels with non-planar menisci. The mathematical equivalence established between thermal resistance and hydrodynamic slip suggests that the generalized solutions derived here may offer new insights into the design of superhydrophobic surfaces with optimized drag-reduction properties. Moreover, the method can be generalised to multiple-phase media (a schematic of three-phase media is illustrated in Figure \ref{fig_3phase}).

\section*{Acknowledgments}
The authors are grateful to Prof.~Darren Crowdy and Prof.~Marc Hodes for useful discussions. EL and AC gratefully acknowledge financial support for the present work provided by {\em QJMAM Fund Awards}. AC is currently supported by JST [Moonshot R\&D][Grant Number JPMJMS25A5].

\bibliographystyle{plain}
\bibliography{thermal}

\end{document}